\documentclass[12pt]{amsart} 

\usepackage{amsthm}

\usepackage{epsf}

\usepackage{amssymb}

\usepackage{latexsym}

\usepackage{amsmath}

  \newcommand{\const}{\rm const}
  
  \newcommand{\supp}{\rm supp}
  
  \newcommand{\Dom}{\rm  Dom}
 \newcommand{\Type}{\rm  Type}

\begin{document}

   \begin{center}

  \   {\bf   Maximal and other operators  in exponential }\\

\vspace{4mm}

   \  {\bf Orlicz and Grand Lebesgue Spaces. }\\

\vspace{4mm}

  \  {\bf Ostrovsky E., Sirota L.}\\

\vspace{4mm}

  Bar-Ilan University, department  of Mathematic and Statistics, ISRAEL, 59200. \\

\vspace{4mm}

E-mails: eugostrovsky@list.ru, \hspace{5mm} sirota3@bezeqint.net \\

\vspace{5mm}

  {\bf Abstract} \\

\vspace{4mm}

 \end{center}

\ We derive in this  preprint the exact up to multiplicative constant non-asymptotical
estimates for  the norms of some non-linear in general case operators, for example, the so-called maximal functional operators,
 in two probabilistic rearrangement invariant norm: exponential Orlicz and Grand Lebesgue Spaces. \par

 \ We will use also the theory of the so-called Grand Lebesgue Spaces (GLS) of measurable  functions. \par

\vspace{5mm}

 \par

\vspace{4mm}


\vspace{5mm}

\section{1. \ Definitions.  Notations. Previous results.  Statement of problem.}

\vspace{4mm}

 \ Let $ \   (X, {B},  {\bf \mu})   \ $ be a probability space: $ \ \mu(X) = 1 . \ $ We will denote by $ \  |f|_p  = |f|L(p) \ $ the ordinary
Lebesgue - Riesz $ \  L(p) \  $ norm of arbitrary measurable numerical valued function $ \   f: X \to R: \  $

$$
|f|_p = |f| L(p) = |f| L_p(X, \mu) := \left[ \int_X |f(x)|^p \ \mu(dx)  \right]^{1/p}, \  p \in [1, \infty).
$$

\vspace{4mm}

 \  {\bf   Definition  1.1.  }  The operator  $ \ Q: L(p) \to L(p), \ p \in (1, \infty), \ $ not necessary to be linear,  acting from
any $ \ L(p) \ $ to one, is said to be of a type $ \ \lambda, \nu; \lambda, \nu = \const, \ \lambda \ge \nu \ge 0, \ $  write

$$
Q \in \Type(\lambda,\nu) = \Type[L(p) \to L(p)](\lambda,\nu), \eqno(1.0)
$$
iff  for some finite constant $ \  Z = Z[Q] \  $ and for certain interval $ \ p \in [1,b), \ b = \const \in (1,\infty) \ $

$$
|Q[f]|_p \le Z \ \frac{p^{\lambda}}{ (p-1)^{\nu}} \ |f|_p, \eqno(1.1)
$$
or equvalently

$$
||Q||[L(p) \to L(p)] \le  Z[Q] \ \frac{p^{\lambda}}{ (p-1)^{\nu}}.  \eqno(1.1a)
$$

\vspace{4mm}

 \ Note that the function $ \  f(\cdot) \ $ in the left-hand side of inequality (1.1)  may be vector - function; moreover,
one can consider  the relation of the form

$$
g(x) = Q(\vec{f}) = Q(f_1, f_2, \ \ldots, f_N), \  N = 1,2,\ldots, \infty,
$$
such that

$$
|g|_p \le  Z(\vec{f}) \cdot \frac{p^{\lambda}}{ (p-1)^{\nu}} \cdot \psi(p) \eqno(1.2)
$$
for certain positive continuous function $ \  \psi = \psi(p), \ p \in [1,b). \ $  The concrete form of these function will be clarified
below. For instance,  one can choose the function $ \  \psi(\cdot) \ $  as a {\it natural function} for the family of the r.v. $ \ \{f_i\}: \ $

$$
\psi(p) := \sup_i |f_i|_p,
$$
if it is finite  at last for one value $  \  p = b, \ b \in (1, \infty]. \ $ \par
 \ This approach may be used for instance in the martingale theory,  see [8],  where $ \  \lambda = \nu = 1,  \  $ and

$$
g = \max_{i=1}^N f_i.
$$

\vspace{4mm}

 \  There are many examples  for  such operators satisfying the estimate (1.2) (or (1.1)):  Doob's  inequality for martingales [8], [13];
 singular integral operators of Hardy-Littlewood type [30], [31],  [24],  Fourier integral operators [27],[30], pseudodifferential
operators [33], theory of Sobolev spaces  [1]   etc.  \par

 \ Note that in the last two examples, as well as in  [31],  [24], \   $ \   \lambda = \nu = 1. \  $\par

 \ Especially many examples are delivered to us the theory  of the so-called  maximal operators, see [1], [8], [27] etc.
For instance,  let $ \  \{ \phi_k \}, \ k = 0,1,2, \ldots \  $ be ordinary complete orthonormal trigonometric system on the set
$ \ [-\pi, \pi] \ $ equipped with {\it normalized} Lebesgue measure $  \ d \mu = dx/(2 \pi) \ $ and let $ \ f \in L(p), \ p \in (1, \infty). \ $
 \ Denote  by

$$
s_M(x) = \sum_{k=0}^M (g, \phi_k) \ \phi_k(x)
$$
the partial Fourier sum for $ \ f(\cdot) \ $ and

$$
g(x)  := \sup_{M \ge 1} |s_M(x)|;
$$
then the inequality (1.2) holds true and herewith $  \ \lambda = 4; \ \nu = 3;  $ see [27]. \par

 \ Note that there are several examples  in which $ \ \lambda > 0, \ $  but $ \   \nu = 0; \ $ see e.g. [32],
[28].\par

 \ Let us turn now our  attention on the ergodic theory, see e.g.   [3], [20], [21] etc.  To be
more concrete, suppose $ \   T = [0,1] \ $ with the classical Lebesgue measure $ \ \mu. \ $  Let $ \   f: T \to  R \   $ be certain  measurable
function.  Denote  as ordinary

$$
f^o(t) := \inf \{  y, \ y > 0: \mu \{s: |f(s)| > y\} \le t \}, \ t \in (0,1]
$$
and

$$
f^{oo}(t) := \frac{1}{t} \int_0^t f^o(x) \ dx.
$$

 \  Introduce for  arbitrary rearrangement invariant (r.i.) space $ \ E \ $  builded over $ \ (T, \mu) \ $  by $ \  H(E) \ $  another
(complete) r.i. space  as follows

$$
H(E) :=  \{ f, \ f \in L_1: f^{oo} \in E   \} \ \subset E
$$
equipped with the norm

$$
||f||H(E) \stackrel{def}{=} ||f^{oo}||E. \eqno(1.3)
$$

 \ Further, let an operator $ \  A \ $ be  an $ \  L_1 \ - \ L_{\infty} \  $  contraction, for instance,

$$
A = A_{\theta} =  A_{\theta}[f](t) = f(\theta(t)),
$$
where $ \ \theta(\cdot) \ $  is an invertible ergodic measure preserving transformation of the set $ \ [0,1]. \ $
 \  Define the following {\it maximal} Dunford-Schwartz operator, not necessary to be linear

$$
B_A[f](t) := \sup_{n = 2,3, \ldots} \left[ \frac{1}{n} \ \sum_{k=0}^{n-1} A^k [f](t)    \right]. \eqno(1.4)
$$
 \ M.Braverman in [3]  proved  that

$$
||B_A[f]||H(E) \le ||f||E. \eqno(1.5)
$$
 \ In particular, if $  \  E = L_p(T), \ 1 < p < \infty, \ $ the estimate (1.5) takes the form

$$
|B_A[f]|_p \le  \frac{p}{p-1} \ |f|_p, \ 1 < p \le \infty.  \eqno(1.6)
$$
 \ So, the inequality (1.1)  is satisfied for the operator $ \ Q = B_A \ $ again with the parameters $  \ \lambda = \nu = 1. \ $  \par

 \ The lower bounds  for the inequalities of the form   (1.1), (1.2), i.e. the lower bounds for the operator $ \ Q \ $  with at the same
parameters $ \ \lambda, \nu \ $ may be found, for instance, in [11], [15]. \par

\vspace{4mm}

\ Notice [3], [9] that there are  rearrangement invariant spaces  $ \ E \ $ for which the norms $ \   ||\cdot||E \ $ and
$ \   ||\cdot||H(E) \ $  are not equivalent. For example,

$$
H(L \ln^{n-1} L) =  L \ln^n L, \ n \ge 1,
$$
see  [3], proposition 1.2.\par

\vspace{4mm}

 \ {\bf  We intend in this preprint to extend the inequality (1.1) (or (1.2)) into the wide class of another rearrangement invariant  Banach
functional spaces: exponential Orlicz spaces and into Grand Lebesgue Spaces. } \par

 \vspace{4mm}

 \ In detail, let $ \  Y_1, \ Y_2 \ $ be  two rearrangement invariant (r.i.) Banach  functional spaces over $ \   (X, {B},  {\bf \mu}), \ $
in particular, Orlicz spaces or Grand Lebesgue ones. We set ourselves the goal to estimate of the correspondent operator norms

$$
||Q||[Y_1 \to Y_2] = \sup_{0 \ne f \in Y_1} \left[ \ \frac{||Qf||Y_2}{||f||Y_1} \ \right]. \eqno(1.7)
$$

\vspace{4mm}

 \section{ 2. Grand Lebesgue Spaces (GLS).}

\vspace{4mm}

  \ Let  $ \   (X, {B},  {\bf \mu})   \ $ be again the source probability space.  Let
also $  \psi = \psi(p), \ p \in [1,b), \ b = \const \in (1, \infty] $ be certain bounded
from below: $ \inf \psi(p) > 0 $ continuous inside the semi - open interval  $ p \in [1,b) $
numerical valued function. We can and will suppose  without loss of generality

 $$
\ \inf_{p \in [1,b)} \psi(p) = 1 \eqno(2.0)
$$
and $  \ b = \sup \{ p, \ \psi(p) < \infty  \}, $
 so that   $ \supp \ \psi = [1, b) $ or  $ \supp \ \psi = [1, b]. $ The set of all such a functions will be
denoted by  $ \ \Psi(b) = \{ \psi(\cdot)  \}; \ \Psi := \Psi(\infty).  $ \par

\vspace{4mm}

 \ By definition, the (Banach) Grand Lebesgue Space (GLS)    $  \ G \psi  = G\psi(b) $
consists on all the real (or complex) numerical valued measurable functions
(random variables, r.v.)   $   \  f: X \to R \ $  defined on our probability space and having a finite norm

$$
|| \ f \ || = ||f||G\psi \stackrel{def}{=} \sup_{p \in [1,b)} \left[ \frac{|f|_p}{\psi(p)} \right]. \eqno(2.1)
$$
 \ The function $ \  \psi = \psi(p) \  $ is said to be {\it  generating function } for this space. \par

\  Furthermore,  let now $  \eta = \eta(z), \ z \in S $ be arbitrary family  of random variables  defined on any set $ \ z \in S \ $ such that

$$
\exists b = \const\in (1,\infty], \ \forall p \in [1,b)  \ \Rightarrow  \psi_S(p) := \sup_{z \in S} |\eta(z)|_p  < \infty.
$$
 \ The function $  p \to \psi_S(p)  $ is named as a {\it  natural} function for the  family  of random variables $  S.  $  Obviously,

$$
\sup_{z \in S} ||\eta(z)||G\Psi_S = 1.
$$

 \ The family $ \ S \ $ may consists on the unique r.v., say $  \  \Delta: \ $

$$
\psi_{\Delta}(p):= |\Delta|_p,
$$
if of course  the last function is finite for some value $ \  p = p_0 > 1. \  $\par
  \ Note that the last condition is satisfied if for instance the r.v. $ \  \Delta \ $ satisfies the so-called Kramer's
condition; the inverse proposition is not true. \par
 \ The generating $ \ \psi(\cdot) \ $ function in (1.2)  may be introduced for instance as natural one for some famoly of a functions. \par

 \ These spaces are Banach functional space, are complete, and rearrangement
invariant in the classical sense, see [2], chapters 1, 2; and were investigated in
particular in many works, see e.g.
 [4], [5], [6], [14], [16], [18], [22], chapters 1,2; [23], [24]  etc. We refer
here some used in the sequel facts about these spaces and supplement more. \par

 \ The so-called tail function $ \ T_{f}(y), \ y \ge 0 \ $ for arbitrary (measurable)  numerical valued function
 $ \  f \ $ is defined as usually

$$
T_{f}(y) \stackrel{def}{=}  \max ( {\bf \mu}(f \ge y), \  {\bf \mu}(f \le -y) ), \ y \ge 0.
$$

 \ It is known that

$$
|f|^p_p= \int_X |f|^p(x) \ \mu(dx) = p \int_0^{\infty} y^{p-1} \ T_f(y) \ dy
$$

and if  $  \ f \in G\psi, \ f  \ne 0, $ then

$$
T_{f}(y) \le  \exp \left( -v_{\psi}^*(\ln(y/||f||G\psi)   \right),  \  y \ge e \ ||f ||G\psi, \eqno(2.2)
$$
where

$$
v(p) = v_{\psi}(p) := p \ \ln \psi(p).
$$

 \ Here and in the sequel the operator  (non - linear) $ \ f \to f^* \ $   will denote the famous Young-Fenchel transform

$$
f^*(u) \stackrel{def}{=} \sup_{x \in \Dom(f)} (x u - f(x)).
$$

 \ Conversely, the last inequality may be reversed in the following version: if

$$
T_{\zeta}(y)  \le  \exp \left( - v_{\psi}^* (\ln (u/K) \right), \ u \ge e \ K,
$$
and if the auxiliary function  $  \ v(p) = v_{\psi}(p)   $ is positive, finite for all the values $ \ p \in [1, \infty),  $ continuous,
convex and such that

$$
\lim_{p \to \infty} \psi(p) = \infty,
$$
then $  \ \zeta \in G(\psi) \  $ and besides $  \ ||\zeta|| \le C(\psi) \cdot K.  $\par

\vspace{4mm}

 \ Let us consider the so-called {\it exponential} Orlicz space $  L(M) $ builded over
source probability space with correspondent Young-Orlicz function

$$
M(y) = M[\psi](y) = \exp \left(  v_{\psi}^*(\ln |y|), \right) \ |y| \ge e;  \ M(y) = C y^2, \ |y| < e.
$$

 \ The exponentiality implies in particular that the Orlicz space  $ \ L(M) \ $ is not separable  as long as the correspondent
Young-Orlicz function  $  \ M(y) = M[\psi](y) \ $ does not satisfy the $ \ \Delta_2 \ $ condition. \par

 \ The Orlicz $ \ ||\cdot||L(M) = ||\cdot||L(M[\psi](\cdot)) \  $ and $ \ ||\cdot|| G\psi \ $ norms are quite equivalent:

$$
||f||G\psi  \le C_1 ||f||L(M) \le C_2 ||f||G\psi,
$$

$$
 0 < C_1 = C_1(\psi) < C_2 = C_2(\psi) < \infty. \eqno(2.3)
$$

\vspace{4mm}

 \  Furthermore,  let now $  \eta = \eta(z), \ z \in W $ be arbitrary family  of measurable functions (random variables)  defined on
any set $  \ W  \ $ such that

$$
\exists b = \const\in (1,\infty], \ \forall p \in [1,b)  \ \Rightarrow  \psi_W(p) := \sup_{z \in W} |\eta(z)|_p  < \infty. \eqno(2.4)
$$
 \ The function $  p \to \psi_W(p)  $ is named as a {\it  natural} function for the  family  of random variables $  W.  $  Obviously,

$$
\sup_{z \in W} ||\eta(z)||G\Psi_W = 1.
$$

 \ The family $ \ W \ $ may consists on the unique r.v., say $  \  \Delta: \ $

$$
\psi_{\Delta}(p):= |\Delta|_p,  \eqno(2.5)
$$
if of course  the last function is finite for some value $ \  p = p_0 > 1. \  $\par
  \ Note that the last condition is satisfied if for instance the r.v. $ \  \zeta \ $ satisfies the so-called Kramer's
condition; the inverse proposition is not true. \par

\vspace{4mm}

 \ {\bf Example 2.0.}  \ Let us consider also the so - called {\it degenerate} $  \ \Psi \ - \ $ function $ \ \psi_{(r)}(p), \ $ where $  r = \const \in [1,\infty): $

$$
\psi_{(r)}(p) \stackrel{def}{=} 1,  \ p \in [1,r];
$$
so that the corresponent value $  b = b(r) $  is equal to $  r. $  One can  extrapolate formally this function onto the whole  semi-axis $  R^1_+: $

$$
\psi_{(r)}(p)  := \infty, \ p > r.
$$

 \  The classical Lebesgue-Riesz $ L_r  $ norm for the r.v. $  \eta $ is  quite equal to the GLS norm $  ||\eta|| G\psi_{(r)}: $

$$
|\eta|_r = ||\eta|| G\psi_{(r)}.
$$
 \ Thus, the ordinary Lebesgue-Riesz spaces are particular, more precisely, extremal cases of the Grand-Lebesgue ones. \par

\vspace{4mm}

  \ {\bf Example 2.1.}  For instance, let  $ \psi $  function has a form

$$
\psi(p) = \psi_m(p) = p^{1/m}, \ m = \const > 0.  \eqno(2.6)
$$

 \ The function $ \ f: X \to R \ $ belongs to the space $ \ G\psi_m: \ $

$$
||f||G\psi_m = \sup_{p \ge 1} \left\{ \  \frac{|f|_p}{p^{1/m}}  \ \right\} < \infty
$$
if and only if  the correspondent tail estimate is follow:

$$
\exists  V = V(m) > 0 \ \Rightarrow \   T_f(y)  \le \exp \left\{ -  (y/V(m))^m   \right\}, \ y \ge 0.
$$

 \  The correspondent  Young-Orlicz function for the space $ \ G\psi_m \ $  has a form

$$
M_m(y) = \exp \left(  |y|^m  \right), \ |y| > 1; \ M_m(y) = e \ y^2, \ |y| \le 1.
$$

 \ There holds for arbitrary  function $ \  f \ $

$$
|| f ||G\psi_m  \asymp ||f||L(M_m) \asymp V(m),
$$
if of course as a capasity of the value $ \ V = V(m) \ $ we understand its minimal positive value from the relation (2.7). \par

 \ The case $ \  m = 2  \ $ correspondent to the so-called subgaussian case, i.e. when

$$
T_f(y) \le \exp \left\{ -  (y/V(2))^2   \right\}, \ y > 0.   \eqno(2.7)
$$

\vspace{4mm}

 \ It is presumes as a rule  in addition that the  function $ \  f(\cdot) \ $ has a mean zero:  $ \  \int_X f(x) \ \mu(dx) = 0. \ $
More examples may be found in [4], [16], [22]. \par

\vspace{4mm}

 \  We   bring  a more general example, see [17].  Let $  \  m = \const > 1 \  $ and define $  \  q = m' = m/( m-1). \  $
Let also $  \  L = L(y), \ y > 0 \ $ be positive  continuous differentiable {\it  slowly varying  }  at infinity function such that

$$
\lim_{\lambda \to \infty} \frac{L(y/L(y))}{L(y)} = 1. \eqno(2.8)
$$
 \ Introduce a following $ \ \psi \ - \ $ function

$$
\psi_{m,L} (p) \stackrel{def}{=} p^{1/m} L^{-1/(m-1)} \left(  p^{ (m-1)^2/m  }  \right\}, \ p \ge 1, \eqno(2.9a)
$$
and a correspondent exponential tail function

$$
T^{(m,L)}(y) \stackrel{def}{=} \exp \left\{  - q^{-1} \ y^q \ L^{ -(q-1)  } \left(y^{q-1} \right)      \right\}, \ y > 0. \eqno(2.9b)
$$

 \ The following implication holds true:

$$
 0 \ne f \in G\psi_{m,L} \ \Longleftrightarrow \exists C = \const \in (0,\infty), \ T_f(y) \le  T^{(m,L)}(y/C). \eqno(2.10)
$$
 \ A particular cases: $ \   L(y) = \ln^r (y+e), \ r = \const, \ y \ge 0; $ then the correspondent generating  functions have a form

$$
\psi_{m,r}(p) = m^{-1} \ p^m  \ \ln^{-r/(m-1)}(p+1),  \eqno(2.11a)
$$
and correspondingly the tail function
$$
T^{m,r}(y) = \exp \left\{ \ - y^q \ (\ln y)^{-(q-1)r}   \ \right\}, \ y \ge e. \eqno(2.11b)
$$

\vspace{4mm}

 \ {\bf Example 2.2.}   Bounded support of generating function. \par

\vspace{4mm}

 \ Introduce the following tail function

$$
T^{<b,\gamma, L>}(x)  \stackrel{def}{=} x^{-b} \ (\ln x)^{\gamma} \ L(\ln x), \ x \ge e, \eqno(2.12)
$$
where as before $ \   L = L(x), \ x \ge 1 \ $ is positive continuous slowly varying function as $ \ x \to \infty, \ $ and

$$
b = \const \in (1, \infty), \ \gamma = \const > -1.
$$

 \ Introduce also the following (correspondent!) $ \ \Psi(b) \ $ function

$$
\psi^{<b,\gamma, L>}(p) \stackrel{def}{=} C_1(b,\gamma,L) \ (b-p)^{ -(\gamma + 1)/b } \ L^{1/b} \left(  \frac{1}{b-p} \right), \
1 \le p < b. \eqno(2.13)
$$

  \ Let  the measurable function $ \ f(\cdot) \ $ be such that

$$
T_f(y) \le T^{<b,\gamma, L>}(y),  \ y \ge e,
$$
then

$$
|f|_p \le  C_2(b,\gamma,L) \ \psi^{<b,\gamma, L>}(p),  \ p \in [1,b)  \eqno(2.14)
$$
or equivalently

$$
 ||f|| \in  G \psi^{<b,\gamma, L>} \ \Longleftrightarrow \  ||f||G \psi^{<b,\gamma, L>}  < \infty. \eqno(2.15)
$$

 \ Conversely, if the estimate (2.14) holds true, then

$$
T_f(y) \le   C_3(b,\gamma,L)  \ y^{-b} \ (\ln y)^{\gamma + 1} \ L(\ln y), \  y \ge e  \eqno(2.16)
$$
or equally

$$
T_f(y) \le T^{<b,\gamma + 1, L>}(y/C_4),  \ y \ge C_4 \ e. \eqno(2.16a)
$$

 \  Notice that there is a  logarithmic ``gap''   as $ \  y \to \infty \ $  between the estimations  (2.15) and (2.16).
 Wherein all the estimates (2.14) and (2.16) are non - improvable, see [17], [18], [24]. \par

\vspace{4mm}

\ {\bf Remark  2.1.} \ These GLS spaces are used for obtaining of an  {\it exponential estimates} for sums of independent random
variables and fields, estimations for non-linear functionals from random fields, theory of Fourier series and transform, theory of operators
  etc., see e.g. [4], [14], [18],  [22],  sections 1.6, 2.1 - 2.5. \par

 \vspace{4mm}

\section{ 3. Main result. The case of equal powers.}

 \vspace{4mm}

 \ We consider in this section  the case when in the relations (1.1) - (1.2)   $ \ \nu = \lambda = \const > 0, \ $ i.e.

$$
|g|_p \le  Z_{\lambda}(\vec{f}) \cdot \frac{p^{\lambda}}{ (p-1)^{\lambda}} \cdot \psi(p) =
 Z \cdot \frac{p^{\lambda}}{ (p-1)^{\lambda}} \cdot \psi(p), \eqno(3.1)
$$
 \ This relation  holds true if  for example  in the relation (1.1) $  \  f \in \  G\psi; \ $  one can assume without loss of generality
 for simplicity $ \ ||f||G\psi = 1, \ Z = 1, \ $ so that $ \  |f|_p \le \psi(p), \ 1 \le p < b. \  $\par

 \vspace{4mm}

 \ {\it Let us introduce some auxiliary constructions.}  Let the function $  \   \psi = \psi(p), \ \psi \in \Psi(b), \ b = \const \in (1, \infty] \  $
be a given. Let also $ \ q \ $ be some fixed number  inside the set $  \  (1, b): \ 1 < q < b. \  $ Suppose the  (measurable)  function
$ \  g = g(x) \ $  satisfies the  inequality (3.1).  We apply the Lyapunov's inequality:   $ \  p \in [1,q] \Rightarrow |g|_p \le |g|_q, \ $
 hence

$$
p \in [1,q] \Rightarrow |g|_p \le |g|_q \le \left[ \frac{q}{q-1} \right]^{\lambda}   \cdot \psi(q). \eqno(3.2a)
$$

 \ We retain the value of  the function $  \  \psi(\cdot) \ $ on the additional set:

$$
p \in (q,b) \Rightarrow |g|_p \le  \left[ \frac{p}{p-1} \right]^{\lambda} \cdot \psi(p). \eqno(3.2b)
$$

 \ Let us introduce the following $ \  \psi \ - \ $  function $ \ \tilde{\psi}(p) = $

$$
\tilde{\psi}_{q, \lambda}(p) := \left[\frac{q}{q-1} \right]^{\lambda}   \cdot \psi(q) \ I( p \in [1,q] ) +
\left[\frac{p}{p-1} \right]^{\lambda} \cdot \psi(p) \ I(p \in (q,b)),  \eqno(3.3)
$$
so that

$$
|g|_p \le Z \cdot \tilde{\psi}_{q, \lambda}(p), \ 1 \le p < b. \eqno(3.4)
$$

 \ Here  and further $ \  I(p \in A)  \ $ denotes the indicator function of the set $ A. $ \par
 \ So, we have eliminated the possible singularity at the point $ \ p \to  1+0. \ $ \par
 \ Let us prove now that

$$
Z^{-1} \ Q(p,q) := \sup_{p \in [1,b)} \frac{\tilde{\psi}_q(p)}{\psi(p) } =: C(b,q,\lambda) < \infty.
$$
 \  We conclude taking into account the restriction $ \ \psi(p)  \ge 1 \ $

$$
Z^{-1} \ Q(p,q) = \frac{\tilde{\psi}_q(p)}{\psi(p) } \le  \left[\frac{q}{q-1} \right]^{\lambda}   \cdot \psi(q) \ I( p \in [1,q] ) +
\left[\frac{p}{p-1} \right]^{\lambda} \ I(p \in (q,b));
$$

$$
\sup_{p \in [1,b)} Z^{-1} \ Q(p,q) \le \max \left(  \frac{q^{\lambda}}{(q-1)^{\lambda}} \psi(q), \  \left[  \frac{q}{q-1} \right]^{\lambda} \right) =
 \left[\frac{q}{q-1} \right]^{\lambda} \ \psi(q). \eqno(3.5)
$$

\vspace{4mm}

 \ Further, let us denote

$$
K_{\lambda}[\psi,b] := \inf_{q \in (1,b)}  \left\{\frac{q^{\lambda}\ \psi(q)}{(q-1)^{\lambda}} \right\},  \eqno(3.6)
$$
then $ \ K_{\lambda}[\psi,b] \in [1, \infty);  \ $ and we derive the following estimate

$$
 1 \le  \inf_q \sup_p \left\{  \frac{\tilde{\psi}_q(p)}{\psi(p) } \right\} \le K_{\lambda}[\psi,b]  < \infty. \eqno(3.7)
$$

  \ We get  due to proper choice of the parameter $  \ q: \ $

\vspace{4mm}

 \ {\bf  Proposition 3.1.} We propose under formulated above notations and conditions, in particular, condition (3.1)

$$
||g||G\psi \le Z \cdot   K_{\lambda}[\psi,b] < \infty. \eqno(3.8)
$$

\vspace{4mm}

 \ One can give a very simple upper estimate for the value $ \ K_{\lambda}[\psi,b];  \ $ indeed, we choose in (3.6)
$  \ q = 2 \ $ in the case when $ \ b > 2  \ $  and $ \ q = (b+1)/2  \ $ if $ \  b \in (1,2]; \  $ we get

$$
K_{\lambda}[\psi,b] \le \left[ \frac{b+1}{b-1} \right]^{\lambda} \ \psi \left( \frac{b+1}{2} \right) \ I(b \in (1,2]) +
2^{\lambda} \ \psi(2) \ I(b > 2).
$$

 \ As a slight consequence: if $ \ b < \infty, \ $ then

$$
K_{\lambda}[\psi,b] \le \frac{C(b,\lambda,\psi)}{(b-1)^{\lambda}}.
$$
where  $ \ C(b,\lambda,\psi) \ $ is continuous bounded function relative the variable $ \ b \ $  in arbitrary finite segment
$ \ 1 < b \le V, \ V = \const  < \infty. $ \par

\vspace{4mm}

\ {\bf Example 3.a.} Let

$$
\psi(p) = \psi_m(p) = p^{1/m}, \ m = \const > 0, \ p \in [1, \infty).
$$

 \ We obtain after simple calculations

$$
K_m(\lambda) := \inf_{q \in (1, \infty)} \left[\frac{q^{\lambda} \ \psi_m(q)}{(q-1)^{\lambda}} \right]=
\inf_{q \in (1, \infty)} \left[ \frac{q^{\lambda + 1/m}}{(q-1)^{\lambda}} \right] =
$$

$$
m^{1/m} \cdot (\lambda + 1/m)^{\lambda + 1/m} \cdot \lambda^{-\lambda}. \eqno(3.9)
$$
 \  In particular,

$$
 K_m(1) = m^{-1 } \ (m+1)^{1 + 1/m}, \ m > 0.
$$

 \ Note by the way  $ \ \forall \lambda > 0 \ \Rightarrow \lim_{m \to \infty} K_m (\lambda) = 1. \ $ \par

\vspace{4mm}

{\bf Example 3.b.} Let $ \ b = \const > 1; \  \beta = \const > 0. $ Define the  following tail function

$$
T [b,\beta](y) : =   C \ y^{-b} \ (\ln y)^{\beta b - 1}, \ y \ge e,
$$
and the following $ \ \Psi(b) \ $ function with bounded support

$$
\psi[b,\beta](p) =    \left[\frac{b-p}{b-1} \right]^{-\beta}, \ p \in [1,b); \  \psi[b,\beta](p)= \infty, \ p \ge b.
$$

 \ The  tail inequality of the form

$$
T_{\eta}(y) \le T[b,\beta](y), \ y \ge e
$$
entails  the inclusion $ \  \eta \in G\psi[b,\beta].  \ $  The inverse conclusion is not true. \par

 \ We find after come computations

$$
K^{b, \beta}(\lambda)  := \inf_{q \in (1,b)} \left[ \frac{q^{\lambda}}{ (q-1)^{\lambda}} \cdot (b - q)^{-\beta}   \right] \le
\frac{(\lambda b + \beta)^{\lambda} \cdot
( \lambda + \beta)^{\beta}}{\lambda^{\lambda} \ \beta^{\beta} \ (b-1)^{\lambda + \beta}}. \eqno(3.10)
$$

\vspace{4mm}

{\bf Example 3.c.} Let now $  \  \psi(p) = \psi_{(r)}(p), \ r = \const > 1. \  $ It is easily to calculate

$$
K_{(r)}(\lambda) := \inf_{q > 1} \left[ \frac{q^{\lambda} \ \psi_{(r)}(q)}{(q-1)^{\lambda}} \right]=  \left[\frac{r}{r-1} \right]^{\lambda}.
$$

\vspace{4mm}

 \ Let us return to the theory of operators,  see (1.1), (1.2). Namely, assume  the operator $ \  Q \  $ satisfies the inequality
(1.1) or more generally (1.2). It follows immediately from proposition (3.1) the following statement. \par

\vspace{4mm}

\ {\bf Theorem 3.1.}  Suppose the function $ \ f(\cdot) \ $ belongs to the space $ \ G\psi \ $ for some generating function
$ \ \psi \ $ from the set $ \ \Psi(b), \ 1 < b \le \infty. \ $ Our statement:
the  function $ \  g = Q[f] \ $ from the relations (1.1)  (or (1.2))
 belongs to at the same Grand Lebesgue Space $ \ G\psi, \  $ or
equivalently to the correspondent exponential Orlicz space $  \ L \left(M_{\psi} \right):  $

$$
  ||g||G\psi \le   K_{\lambda}[\psi,b] \  Z(Q) \ ||f||G\psi, \eqno(3.11)
$$
or equally in the terms of  exponential Orlicz spaces

$$
  ||g|| L \left(M_{\psi} \right) \le  \ C(\psi) \  K_{\lambda}[\psi,b] \  Z(Q) \ ||f|| L \left(M_{\psi} \right). \eqno(3.11a)
$$

\vspace{4mm}

 \ {\bf Remark 3.1.}  The statement of theorem (3.1) may be reformulated as follows. Under at the same conditions:
$ \  f \in G\psi \ $ etc.

$$
||Q[f]||G\psi \le Z \   K_{\lambda}[\psi,b] \ ||f||G\psi \eqno(3.11c)
$$
or equally

$$
||Q(\cdot)||[G\psi \to G\psi] \le  Z \   K_{\lambda}[\psi,b]. \eqno(3.11d)
$$

\vspace{4mm}

 \ {\bf Remark 3.2.} Note that the considered here Young-Orlicz function  $  \ M_{\psi} (y) \ $ does not satisfy the
$ \  \Delta_2  \ $ condition, in contradiction to the considered ones in the book  [19], section 12. \par

\vspace{4mm}

 \ {\bf Example 3.1.} Suppose the  function  $ \ f(\cdot) \ $   from the estimate (1.1)  belongs to the space
$ \ G\psi_m, \ m = \const > 0: $

$$
\sup_{p \ge 1} \left[ \frac{|f|_p}{p^{1/m}}   \right] < \infty \eqno(3.12a)
$$

 or equivalently
$$
 \exists C_1 > 0 \ \Rightarrow T_f(y)  \le \exp(- C_1 \ y^m), \ y \ge 0. \eqno(3.12b)
$$
Then  there exists a positive  finite constant $  \ C_3 = C_3(m,\lambda) \ $ for which

$$
  T_g(y)  \le \exp(- C_3(m,\lambda) \ y^m), \ y \ge 0, \eqno(3.13a)
$$
 or equivalently

$$
\sup_{p \ge 1} \left[ \frac{|g|_p}{p^{1/m}}   \right] < \infty. \eqno(3.13b)
$$

\vspace{4mm}

 \ More generally, let  $  \  L = L(y), \ y > 0 \ $ be the positive continuous differentiable {\it  slowly varying  }  at infinity
function such that

$$
\lim_{\lambda \to \infty} \frac{L(y/L(y))}{L(y)} = 1,
$$
i.e. as in the example 2.1.  Recall the following notation for $ \ \psi \ - \ $ function

$$
\psi_{m,L} (p) \stackrel{def}{=} p^{1/m} L^{-1/(m-1)} \left(  p^{ (m-1)^2/m  }  \right\}, \ p \ge 1,  \ m > 1,
$$
and  the  correspondent exponential tail function

$$
T^{(m,L)}(y) \stackrel{def}{=} \exp \left\{  - q^{-1} \ y^q \ L^{ -(q-1)  } \left(y^{q-1} \right)  \right\}, \ y > 0,
$$
where $   \  m = \const > 1,  \ q = m/(m-1). \  $\par
 \ Suppose the  function  $ \ f(\cdot) \ $   from the estimate (1.1) belongs to the space
$ \ G\psi_{m,L}, \ m = \const > 0: $

$$
\sup_{p \ge 1} \left[ \frac{|f|_p}{\psi_{m,L}(p)}   \right] < \infty  \eqno(3.14a)
$$

 or equivalently
$$
 \exists C_1 = C_1(m,L) > 0 \ \Rightarrow T_f(y)  \le T^{(m,L)}(y/C_1). \eqno(3.14b)
$$
Then there exists a positive   constant $  \ C_3 = C_3(m,L,\lambda) \ $ for which

$$
  T_g(y)  \le T^{(m,L)}( y/C_3)  \eqno(3.15a)
$$
 or equivalently

$$
\sup_{p \ge 1} \left[ \frac{|g|_p}{\psi_{m,L}(p)}   \right] < \infty.  \eqno(3.15b)
$$

\vspace{4mm}

{\bf Example 3.2.}  The case of bounded support. \par

\vspace{4mm}

 \ This case is more complicated. Recall the following notation for tail function

$$
T^{<b,\gamma, L>}(x)  \stackrel{def}{=} x^{-b} \ (\ln x)^{\gamma} \ L(\ln x), \ x \ge e,
$$
where as before $ \   L = L(x), \ x \ge 1 \ $ is the positive continuous slowly varying function as $ \ x \to \infty, \ $ and let as before

$$
b = \const \in (1, \infty), \ \gamma = \const > -1.
$$

 \ and recall also notation for the following correspondent $ \ \Psi(b) \ $ function

$$
\psi^{<b,\gamma, L>}(p) \stackrel{def}{=} C_1(b,\gamma,L) \ (b-p)^{ -(\gamma + 1)/b } \ L^{1/b} \left(  \frac{1}{b-p} \right), \
1 \le p < b.
$$

  \ Let  the source  (measurable) function $ \ f(\cdot) \ $ be such that

$$
 ||f|| \in  G \psi^{<b,\gamma, L>} \ \Longleftrightarrow \  ||f||G \psi^{<b,\gamma, L>}  < \infty,  \eqno(3.16)
$$
then also $ ||g|| \in  G \psi^{<b,\gamma, L>}  \ $ and moreover

$$
  ||g||G \psi^{<b,\gamma, L>}  \le  K^{b, \gamma}(\lambda) \  ||f||G \psi^{<b,\gamma, L>} . \eqno(3.17)
$$

\vspace{4mm}

 \ But if  we assume  the following tail restriction  on the function $  \  f \  $

$$
T_f(y) \le T^{<b,\gamma, L>}(y),  \ y \ge e,  \eqno(3.18)
$$
then  we conclude only

$$
T_g(y) \le   C_5(b,\gamma,L)  \ y^{-b} \ (\ln y)^{\gamma + 1} \ L(\ln y), \  y \ge e   \eqno(3.19a)
$$
or equally

$$
T_g(y) \le T^{<b,\gamma + 1, L>}(y/C_6),  \ y \ge C_6 \ e.  \eqno(3.19b)
$$

 \ {\it  Open question: \ } what is the ultimate value instead $ \ ``\gamma + 1'' \ $ in the last estimate? \par

 \vspace{4mm}

\section{ 4. Main result. The case of  different  powers.}

 \vspace{4mm}

 \ Let as before some function $ \  \psi = \psi(p), \ p \in [1,b), \ b = \const \in (1, \infty] \  $ from the set $ \ \Psi(b) \ $
be a given. \ Suppose in this section that in the inequalities (1.1) or (1.2)  $ \  f \in G\psi, \ ||f||G\psi  < \infty, \ \lambda > \nu \ge 0, \ $
and denote $ \  \Delta =  \lambda - \nu;\ (\Delta > 0),  \ $

$$
\zeta(p) = \zeta[\psi, \Delta](p) :=  p^{\Delta} \ \psi(p). \eqno(4.0)
$$
 \ Obviously, $ \  \zeta(\cdot) \in \Psi(b). \ $\par

\ Let for beginning $  \  ||f||G\psi = 1, \  $  then $  \  |f|_p \le \psi(p), \  p \in [1,b). \  $
  We  deduce from the inequality  (1.1) taking into account the estimate $ \ |f|_p \le \psi(p), \ p \in [1,b) \ $
alike the foregoing section  denoting $ \  g = Q[f] \  $

$$
|g|_p \le Z \ \left[  \frac{p}{p-1}  \right]^{\nu} \ p^{\Delta} \ \psi(p) =
 Z \ \left[  \frac{p}{p-1}  \right]^{\nu} \ \zeta[\psi, \Delta](p). \eqno(4.1)
$$

 \ It follows immediately from proposition (3.1)  or theorem 3.1

$$
||g||G\zeta[\psi, \Delta] \le Z \ K_{\nu}[\zeta[\psi,\Delta]](\zeta, b). \eqno(4.2)
$$

\vspace{4mm}

 \ We proved in fact the following result. \par

\vspace{4mm}

\ {\bf Theorem 4.1.}  Suppose as above that the function $ \ f(\cdot) \ $ belongs to the space $ \ G\psi \ $ for some generating function
$ \ \psi \ $ from the set $ \ \Psi(b), \ 1 < b \le \infty. \ $  Let in (1.1) $  \ \lambda > \nu \ge 0.   \  $ Our statement:
the  function $ \  g = Q[f] \ $ from the relations (1.1)  belongs to  the  {\it other}  certain Grand Lebesgue Space $ \ G\zeta, \  $ or
equivalently to the correspondent exponential Orlicz space $  \ L \left(M_{\zeta} \right):  $

$$
||g||G\zeta[\psi, \Delta] \le Z \ K_{\nu}(\zeta, b) \ ||f||G\psi, \eqno(4.3)
$$

or equally in the terms of exponential Orlicz spaces

$$
  ||g|| L \left(M_{\zeta} \right) \le  \ C(\psi) \  K_{\nu}[\zeta[\psi,\Delta]](\nu,b) \  Z(Q) \ ||f|| L \left(M_{\psi} \right). \eqno(4.3a)
$$

\vspace{4mm}

{\bf Remark 4.1.}  In the case $ \  b < \infty \ $  the estimate (4.3) may be simplified as follows. As long as  in this case
$ \   p^{\Delta}  \ \psi(p) \  \le b^{\Delta} \ \psi(p), \ $ we conclude that the operator $ \ Q \ $ acts from the space
$ \   G\psi \  $ into {\it at the same space:}

$$
||Q[f]||G\psi \le  b^{\Delta} \ Z \ K_{\nu} \left[b^{\Delta} \psi, b \right] \ ||f||G\psi.  \eqno(4.4)
$$

 \vspace{4mm}

 \section{ 5. Convergence in the Grand Lebesgue and non-separable Orlicz spaces. }

 \vspace{4mm}

  \ Let us  consider here  the {\it sequence} of the form

$$
g_n(x) = Q_n(\vec{f}_n) = Q(f_1, f_2, \ \ldots, f_n), \  n = 1,2,\ldots, \infty,  \eqno(5.0)
$$
such that for some non - negative  constants   $   \ \lambda, \nu; \ \lambda \ge \nu    \  $

$$
\sup_n |g_n|_p \le  Z(\vec{f}) \cdot \frac{p^{\lambda}}{ (p-1)^{\nu}} \cdot \psi(p) \eqno(5.1)
$$
for certain positive continuous function $ \  \psi = \psi(p), \ p \in [1,b), \ b = \const \in (1, \infty] $  from the set $ \ \Psi(b). \ $ \par

 \ It follows from theorem 4.1 that

$$
\sup_n ||g_n||G\zeta < \infty, \eqno(5.2)
$$
where in the case $  \   \lambda = \nu  \ \Rightarrow \zeta(p) = \psi(p). \  $\par

\vspace{4mm}

 \ {\it   We suppose in addition to (5.1) (or following (5.2) ) that the sequence   } $ \   \{  g_n(\cdot)  \} \ $
{\it  converges in all the norms } $ \  L_p(X, \mu), \ p \in [1,b); \  $

$$
\exists g_{\infty}(x)  \ \stackrel{def}{=} \lim_{n \to \infty} g_n(x) \eqno(5.3a)
$$
such that

$$
\forall p \in [1,b) \ \Rightarrow \lim_{n \to \infty} |g_n - g_{\infty}|_p = 0. \eqno(5.3b)
$$

\vspace{4mm}

 \ Our claim in this section is investigation under formulated before condition the problem of convergence
$ \   g_n \to g_{\infty} \ $ in more strong norms, concrete: in the CLS sense or correspondingly in Orlicz  spaces norms. \par

\vspace{4mm}

 \ The simplest example of (5.2)-(5.3a), (5.3b) give us the theory of martingales. It makes sense to dwell on this in more detail.

\vspace{4mm}

 \ This approach may be used for instance in the martingale theory, see [8], where $ \  \lambda = \nu = 1,  \  $
and

$$
g_n = \max_{i=1}^n f_i.
$$
where $ \ \{ f_i  \}  \ $ is a centered martingale (or semi-martingale) sequence relative certain  filtration $ \ \{ F_i \}: \ $ \par

$$
|g_n|_p \le \frac{p}{p-1} \  |f_n|_p, \ p \in (1, b),
$$
if of course the right-hand side is finite. \par

\vspace{4mm}

 \ J. Neveu  proved in  [19], pp. 209-220  that if the Orlicz space $ \  L(M) \ $ builded over our probability space
with correspondent Young-Orlicz function  $ \ M(\cdot) \ $
satisfying the $ \ \Delta_2 \ $ condition,  or equvalently if the space $ \  L(M) \ $ is separable,

$$
 \lim_{t \to 0+} M(t)/t  = 0
$$

and

$$
\sup_{n = 1,2, \ldots} ||f_n||L(M) < \infty,
$$
then there esists almost ewerywhere a limit

$$
\lim_{n \to \infty} f_n =: f_{\infty}
$$
and the convergence in (4.2) take place also in the $ \  L(M) \ $ norm:

$$
\lim_{n \to \infty} || \ f_n -  f_{\infty} \ ||L(M) = 0. \eqno(5.4)
$$

 \ We must  first of all recall some definitions and facts about comparison of GLS from an article  [24].
Let  $ \   \psi, \nu \  $ be two functions from the set $ \  G\psi(b), \ b \in (1, \infty]. \ $ We will write
$ \ \psi << \nu, \ $ or equally $ \  \nu >> \psi, \ $ iff

$$
\lim_{p \to b-0} \frac{\psi(p)}{\nu(p)} = 0, \ b < \infty;\eqno(5.5a)
$$

$$
\lim_{p \to \infty} \frac{\psi(p)}{\nu(p)} = 0, \ b = \infty. \eqno(5.5b)
$$

 \ There exists an equivalent version (and notion) for Young-Orlicz function, see [26], chapters 2,3. \par

\vspace{4mm}

 \ {\bf Theorem 4.1.} Assume the formulated above notations, conditions  (5.3a), (5.3b)   remains true.  Our statement:
for arbitrary $ \  \Psi(b) \ $ function  $ \ \tau = \tau(p), \ 1 \le p < b \ $ such that $  \  \tau << \zeta \ $

$$
\lim_{n \to \infty} || f_n - f_{\infty} \ ||G\tau = 0, \eqno(5.6)
$$
i.e. the sequence $  \ f_n \ $ converges  not only almost surely but also in arbitrary $   \  G\tau \  $ norm
for which $ \ \tau << \zeta. \ $ \par

\vspace{4mm}

{\bf Proof} is very simple. \ The needed convergense $ \ f_n, \ n \to \infty \ $ in the $ \  G\tau  \ $ norm  follows immediately
from one of the main results of the article [24], p. 238. \par

\vspace{4mm}

\section{ 6. Concluding remarks. } \par

 \vspace{4mm}

 \  {\bf A.}  One can consider a more general case as in (1.1), (1.2):

$$
|Q[f]|_p \le W(p) \ |f|_p,  \ p \in [1,b), \ f \in G\psi, \ \eqno(6.0)
$$
or equally $ \ g(x) := Q[f](x), $

$$
|g|_p \le W(p) \ \psi(p),  \ p \in [1,b), \ \eqno(6.1)
$$
where $ \  W = W(p), \ p \in (1,b) \  $ is any measurable function,  not necessary to be continuous or bounded. \par

 \ Indeed,  let as above $ \ q \ $ be arbitrary number from the open interval $ \  (1,b): \ q \in (1,b).  \ $ We have using again
the Lyapunov's inequality

$$
p \in [1,q] \ \Rightarrow |g|_p \le I(p \in [1,q]) \ W(q) \ \psi(q),
$$
following

$$
|g|_p \le I(p \in [1,q]) \ W(q) \ \psi(q) + I(p \in (q,b) ) \  W(p) \ \psi(p).
$$
 \ Thus, if we denote $ \ \upsilon(p) = \upsilon[W,\psi](p) := \ $

$$
 \inf_{q \in (1,b)} \left\{  I(p \in [1,q]) \ W(q) \ \psi(q) + I(p \in (q,b) ) \  W(p) \ \psi(p) \right\}: \eqno(6.2)
$$

\vspace{4mm}

\ {\bf  Proposition 6.1.  \ }

$$
||Q[f]||G \upsilon \le ||f||G\psi. \eqno(6.3)
$$

\vspace{4mm}

\ {\bf B.} \  In the case of  martingales the condition of almost surely convergence follows from the boundedness of its moment:
$  \  \sup_n |f_n|_p < \infty, \ \exists p \ge 1, \  $ by virtue of the famous theorem of J.Doob. \par

\vspace{4mm}

\ {\bf C.} \  Lower bounds for the norm of considered operators. \par

\vspace{4mm}

 \ Assume in addition to the estimates (1.1), or (1.2), (3.1),  that

$$
|f|_p \le |Q[f]|_p, \ p \in [1,b), \  b = \const  \in (1, \infty). \eqno(6.4)
$$

 \ The inequality (6.3) is true for example for every maximal operators, in the J.Doob's inequality for  martingales etc. \par
 \  Suppose that the function $ \  \psi(\cdot) \ $ is a natural function for appropriate function $ \ f_0: X \to R: \ \psi(p) = |f_0|_p, \ $
such that $  \  \forall p < b \ \Rightarrow  \psi(p) < \infty. \  $  Then the  relation (6.4) takes the form

$$
\psi(p) \le  |Q[f_0]|_p, \ p \in [1,b), \ \Longleftrightarrow \ |Q[f_0]|_p \ge \psi(p),
$$
hence

$$
||Q[\cdot]||[G\psi \to G\psi] \ge 1. \eqno(6.5)
$$

\vspace{4mm}

 {\bf References.}

\vspace{4mm}

{\bf 1. R.A. Adams.}  {\it Sobolev Spaces.} Academic Press, New York, 1975.\\
{\bf 2. Bennet C., Sharpley R. } {\it Interpolation of operators.} Orlando, Academic
Press Inc., (1988). \\
{\bf 3. M.Braverman, Ben-Zion Rubstein, A.Veksler.}  {\it Domimated ergodic theorems in rearrangement invariant spaces. \ }
Studia Mathematica, {\bf 128}, (2), 145-157, (1998). \\
{\bf 4. Buldygin V.V., Kozachenko Yu.V.} {\it Metric Characterization of Random
Variables and Random Processes. } 1998, Translations of Mathematics Monograph,
AMS, v.188. \\
{\bf 5. A. Fiorenza.} {\it Duality and reflexivity in grand Lebesgue spaces.} Collect.
Math. 51, (2000), 131-148.\\
{\bf 6.  A. Fiorenza and G.E. Karadzhov.} {\it Grand and small Lebesgue spaces and
their analogs.} Consiglio Nationale Delle Ricerche, Instituto per le Applicazioni
del Calcoto Mauro Picone”, Sezione di Napoli, Rapporto tecnico 272/03, (2005). \\
{\bf 7.  Y.Derriennic.}  {\it On integrability of the supremum of ergodic  ratios.} Ann. Probab., {\bf 1,}  (1973), 338 - 340. \\
{\bf 8.  J.L.Doob.} {\it Stochastic Processes.} Wiley, New York, 1953. \\
{\bf 9. D.Gilat.} {\it The best bound in the $ \ L \log L \ $  inequality of Hardy  and Littlewood  and its  martingale  counterpart.}
 Proc. Amer. Math. Soc., 1992, {\bf 97,}  429 - 436. \\
{\bf 10. D.Gilat.} {\it On the ratio of the expected maximum of a martingale and the $ \ Lp \ - \ $ norm of its last term.}
Israel Journal of Mathematics, October 1988, Volume 63, Issue 3, pp 270–280. \\
{\bf 11. L. Grafakos, A. Montgomery-Smith.} {\it Best constants for uncentered maximal functions.}
Bulletin of the London Mathematical Society, 29, (1997), no. 1, 60–64.
{\bf 12. I.Grama, E.Haeusler. }  {\it Large deviations for martingales. } Stoch. Pr.
Appl., 85, (2000), 279-293. \\
{\bf 13.  P.Hall, C.C.Heyde.} {\it  Martingale Limit Theory and Its Applications. } Academic Press. 1980,
 New York, London, Toronto, Sydney.\\
{\bf 14. T. Iwaniec and C. Sbordone.} {\it On the integrability of the Jacobian under
minimal hypotheses.} Arch. Rat.Mech. Anal., 119, (1992), 129-143.\\
{\bf 15.  Paata Ivanisvili, Benjamin Jave, and Fedor Nazaro. } {\it  Lower bounds for uncentered maximal functions in
any dimension. } arXiv:1602.05895v1  [math.AP]  18 Feb 2016.
{\bf 16. Kozachenko Yu. V., Ostrovsky E.I.}  (1985). {\it The Banach Spaces of
random Variables of subgaussian Type.} Theory of Probab. and Math. Stat. (in
Russian). Kiev, KSU, 32, 43-57. \\
{\bf 17.  Kozachenko Yu.V.,   Ostrovsky E.,   Sirota L} {\it Relations between exponential tails, moments and
moment generating functions for random variables and vectors.} \\
arXiv:1701.01901v1  [math.FA]  8 Jan 2017 \\
{\bf 18. E. Liflyand, E. Ostrovsky and L. Sirota.} {\it Structural properties of Bilateral Grand Lebesque Spaces. }
Turk. Journal of Math., 34, (2010), 207-219. TUBITAK, doi:10.3906/mat-0812-8 \\
{\bf 19. Neveu J.} {\it Discrete - parameter Martingales.} North - Holland, Amsterdam, 1975.\\
{\bf 20.  D.S.Ornstein. } {\it  A remark on the Birkhoff ergodic theorem.}  Illinois  J. Math., {\,15}  (1971), 77 - 79. \\
{\bf 21.  Adam Ose'kowski.} {\it  Sharp logarithmic inequalities for two Hardy-type operators. }
Zeitschrift f\''ur Analysis und ihre Anwendungen. European Mathematical Society,  Journal for Analysis and its Applications, Volume 20, (2012),
Monografie Matematyczne 72 (2012), Birkhäuser, 462 pp.  DOI: 10.4171/ZAA/XXX. \\
{\bf 22. Ostrovsky E.I.} (1999). {\it Exponential estimations for Random Fields and
its applications,} (in Russian). Moscow-Obninsk, OINPE.\\
{\bf 23. Ostrovsky E. and Sirota L.}  {\it Sharp moment estimates for polynomial martingales.}\\
arXiv:1410.0739v1 [math.PR] 3 Oct 2014 \\
{\bf 24. Ostrovsky E. and Sirota L.}  {\it  Moment Banach spaces: theory and applications. }
HAIT Journal of Science and Engineering C, Volume 4, Issues 1-2, pp. 233-262.\\
{\bf 25. De La Pena V.H. }  {\it A general class of exponential inequalities for martingales and ratios.} 1999,
Ann. Probab., 36, 1902-1938. \\
{\bf 26. Rao M.M.,  Ren Z.D.}  {\it Theory of Orlicz spaces.} (1991), Pure and Applied Mathematics, Marcel Dekker,
ISBN 0-8247-8478-2. \\
{\bf  27.  Juan Arias de Reyna.}  {\it Pointwise Convergence of Fourier Series.}  Springer, 2002. \\
{\bf 28.  Derek W. Robinson.} {\it Hardy and Rellich inequalities on the complement of convex sets.}
arXiv:1704.03625v1  [math.AP]  12 Apr 2017 \\
{\bf 29.  Secchi,  S.,  Smets,  D., and Willem,  M.} {\it  Remarks on a Hardy–Sobolev inequality.}
C. R. Acad. Sci Paris, Serie I, 336, (2003), 811–815.\\
{\bf 30.   E M. Stein and G.Weiss..} {\it  Introduction to Fourier Analysis  on Euclidean Spaces.}, (1971),
Princeton Univ. Press,   Princeton, N.J.  \\
{\bf 31.  Elias M. Stein.} {\it Singular Integrals and Differentiability Properties of Functions.}  (PMS-30), (1978), Oxford. \\
{\bf  32. Elias M. Stein.} {\it Maximal functions: spherical means.} Proceedings of the National Academy of Sciences,
73 (7),  (1976), 2174–2175.\\
{\bf 33.  Michael E. Taylor.} {\it Pseudodifferential Operators  And Nonlinear PDE.} Springer, 1970, Published by: Princeton
University Press. \\

\vspace{4mm}

\end{document}